\documentclass[10pt]{amsart}
\usepackage{latexsym}
\usepackage{amsfonts}

\usepackage{amsmath,amsthm,amssymb}

\title[Densities in free groups and $\mathbb{Z}^k$] {Densities in free
  groups and $\mathbb{Z}^k$,\\ Visible Points and Test Elements}

\author[I.~Kapovich]{Ilya Kapovich}

\address{\tt Department of Mathematics, University of Illinois at
  Urbana-Champaign, 1409 West Green Street, Urbana, IL 61801, USA
  \newline http://www.math.uiuc.edu/\~{}kapovich/} \email{\tt
  kapovich@math.uiuc.edu}

\author[I.~Rivin]{Igor Rivin} \address{\tt Department of Mathematics,
  Temple University, Philadelphia, PA 19122, USA} \email{\tt
  rivin@euclid.math.temple.edu} \curraddr{Mathematics Department,
  Princeton University, Fine Hall, Washington Rd, Princeton NJ 08544,
  U.S.A.}

\author[P.~Schupp]{Paul Schupp}

\address{\tt Department of Mathematics, University of Illinois at
  Urbana-Champaign, 1409 West Green Street, Urbana, IL 61801, USA
  \newline http://www.math.uiuc.edu/People/schupp.html}
\email{schupp@math.uiuc.edu}

\author[V.~Shpilrain]{Vladimir Shpilrain} \address{\tt Department of
  Mathematics, The City College of New York, New York, NY 10031, USA
  \newline http://www.sci.ccny.cuny.edu/\~{}shpil} \email{\tt
  shpil@groups.sci.ccny.cuny.edu}

\newtheorem{theor}{Theorem}

\newtheorem{thm}{Theorem}[section] \newtheorem{lem}[thm]{Lemma}
 
\newtheorem{prop}[thm]{Proposition} \theoremstyle{definition}
\newtheorem{defn}[thm]{Definition}
\newtheorem{notation}[thm]{Notation}
\newtheorem{conv}[thm]{Convention} \newtheorem{rem}[thm]{Remark}
\newtheorem{exmp}[thm]{Example}

 \newtheorem{prob}[thm]{Problem}

\begin{document}

\begin{abstract}
  In this article we relate two different densities.  Let $F_k$ be the
  free group of finite rank $k \ge 2$ and let $\alpha$ be the
  abelianization map from $F_k$ onto $ \mathbb{Z}^k$.  We prove that
  if $S \subseteq \mathbb{Z}^k$ is invariant under the natural action
  of $SL(k, \mathbb{Z})$ then the asymptotic density of $S$ in
  $\mathbb Z^k$ and the annular density of its full preimage
  $\alpha^{-1}(S)$ in $F_k$ are equal. This implies, in particular,
  that for every integer $t\ge 1$, the annular density of the set of
  elements in $F_k$ that map to $t$-th powers of primitive elements in
  $\mathbb{Z}^k$ is equal to to $\frac{1}{t^k\zeta(k)}$, where $\zeta$
  is the Riemann zeta-function.

  An element $g$ of a group $G$ is called a \emph{test element} if
  every endomorphism of $G$ which fixes $g$ is an automorphism of $G$.
  As an application of the result above we prove that the annular
  density of the set of all test elements in the free group $F(a,b)$
  of rank two is $1-\frac{6}{\pi^2}$.  Equivalently, this shows that
  the union of all proper retracts in $F(a,b)$ has annular density
  $\frac{6}{\pi^2}$.  Thus being a test element in $F(a,b)$ is an
  ``intermediate property'' in the sense that the probability of being
  a test element is strictly between $0$ and $1$.
\end{abstract}


\maketitle

\section{Introduction}\label{intro}

The idea of genericity and generic-case behavior in finitely presented
groups was introduced by Gromov~\cite{Grom,Grom1} and is currently the
subject of much research. (See, for
example~\cite{AO,A1,A2,BMS,Ch95,Ch00,Che98,Grom2,KMSS,KMSS1,KSS,KS,KSn,Oliv,Ol92,Z}.)
Looking at the properties of random groups led Gromov~\cite{Grom2} to
a probabilistic proof that there exists a finitely presented group
that is not uniformly embeddable in a Hilbert space. It also turns out
that random group-theoretic objects exhibit various kinds of algebraic
rigidity properties.  In particular, Kapovich, Schupp and
Shpilrain~\cite{KSS} proved that a random cyclically reduced element
of a free group $F=F(A)$ is of minimal length in its $Aut(F)$-orbit
and that such an element has a trivial stabilizer in $Out(F)$.
Moreover, it turns out~\cite{KSS} that random one-relator groups
satisfy a strong Mostow-type rigidity.  Specifically, two random
one-relator groups $G_r=\langle a_1,\dots ,a_k| r \rangle$ and $G_s
=\langle a_1,\dots ,a_k| s \rangle$ are isomorphic if and only if
their Cayley graphs on the \emph{given} set of generators
$\{a_1,\dots, a_k\}$ are isomorphic as labelled graphs where the graph
isomorphism is only allowed to permute the label set $\{a_1,\dots,
a_k\}^{\pm 1}$.

The most straightforward definition of ``genericity'' is based on the
notion of ``asymptotic density''.

\begin{defn}[Asymptotic density]\label{defn:as}
  Suppose that $T$ is a countable set and that $\ell:T\to \mathbb N$
  is a function (referred to as \emph{length}) such that for every
  $n\in \mathbb N$ the set $\{x\in T: \ell(x)\le n\}$ is finite.  If
  $X\subseteq T$ and $n\ge 0$, we denote $\rho_{\ell}(n,X):=\#\{x\in
  X: \ell(x)\le n\}$ and $\gamma_\ell(n,S)=\#\{x\in X: \ell(x)=n\}$.
  
  Let $S\subseteq T$. The \emph{asymptotic density} of $S$ in $T$ is
\[
\overline{\rho}_{T,\ell}(S):=\limsup_{n\to\infty}\frac{\#\{x\in S:
  \ell(x)\le n\}}{\#\{x\in T: \ell(x)\le n\}}=\\
\limsup_{n\to\infty}\frac{\rho_\ell(n,S)}{\rho_\ell(n,T)},
\]
where we treat a fraction $\frac{0}{0}$, if it occurs, as $0$.

 If the actual
limit exists, we denote it by $\rho_{T,\ell}(S)$ and call this
limit the \emph{strict asymptotic density} of $S$ in $T$. We say
that $S$ is \emph{generic in $T$ with respect to $\ell$} if
$\rho_{T,\ell}(S)=1$ and that $S$ is \emph{negligible} in $T$ if
$\rho_{T,\ell}(S)=0$.
\end{defn}

If $S$ is $T$-generic then the probability that a uniformly chosen
element of $T$ of length at most $n$ belongs to $S$ tends to $1$ as
$n$ tends to infinity.

It turns out that a different density, recording the proportions of a
set in two successive spheres, is sometimes more suitable for subsets
of a free group.

\begin{defn}[Annular Density]
  Let $T,S,\ell$ be as in Definition~\ref{defn:as}. The \emph{annular
    density} of $S$ in $T$ with respect to $\ell$ is:
\begin{gather*}
  \overline\sigma_{T,\ell}(S):=
  \limsup_{n\to\infty}\frac{1}{2}\big(\frac{\#\{x\in S: \ell(x)=
    n-1\}}{\#\{x\in T: \ell(x)= n-1\}}+\frac{\#\{x\in S: \ell(x)=
    n\}}{\#\{x\in T: \ell(x)= n\}}\big)=\\  
  \limsup_{n\to\infty}\frac{1}{2}\big(\frac{\gamma_\ell(n-1,S)}{\gamma_\ell(n-1,T)}+
  \frac{\gamma_\ell(n,S)}{\gamma_\ell(n,T)}\big),
\end{gather*}
where we treat a fraction $\frac{0}{0}$, if it occurs, as $0$.  Again,
if the actual limit exists, we denote this limit by
$\sigma_{T,\ell}(S)$ and call it the \emph{strict annular density} of
$S$ in $T$ with respect to $\ell$.
\end{defn}

\begin{conv}
  Throughout this paper $F=F(A)$ will be a free group of rank $k \ge
  2$ with a fixed finite basis $A=\{a_1,\dots,a_k\}$.  If $w\in F$
  then $|w|_A$ denotes the freely reduced length of $w$ with respect
  to the basis $A$. In discussing the density (asymptotic or annular)
  of subsets of $F$ using the notation above, we will assume that the
  ambient set is $T=F$ and that the length function $\ell(w)$ is
  $|w|_A$.  If $S\subseteq F$ we denote its asymptotic and annular
  densities by $\overline{\rho}_A(S)$ and $\overline{\sigma}_A(S)$
  respectively, and if the strict asymptotic density or the strict
  annular density exist we denote them by $\rho_{A}(S)$ and
  $\sigma_A(S)$ respectively.  Also, we denote
  $\gamma_A(S):=\gamma_{\ell}(S)$ and $\rho_A(n,S):=\rho_\ell(n,S)$ in
  this case.

  For subsets of $\mathbb Z^k$ a length function $\ell: \mathbb Z^k\to\mathbb R$
  will usually be the restriction to $\mathbb Z^k$ of
  $||.||_p$-norm from $\mathbb R^k$ for some $1\le p\le \infty$.
  In this case for $S\subseteq Z^k$ we denote the corresponding asymptotic density of $S$
  in $T=\mathbb Z^k$ by $\overline{\rho}_p(S)$ and if
  the strict asymptotic density exists, we denote it by
  $\rho_{p}(S)$.
\end{conv}

It is not hard to see that if for a subset $S\subseteq F$ the strict
asymptotic density $\rho_A(S)$ exists then the strict annular density
$\sigma_A(S)$ also exists and in fact $\sigma_A(S)=\rho_A(S)$. Namely,
since the sizes of both the balls and the spheres in $F(A)$ grow as
constant multiples of $(2k-1)^n$, if the strict asymptotic density
$\rho_A(S)$ exists, then the limit $\lim_{n\to\infty}
\frac{\gamma_A(n,S)}{\gamma_A(n,F)}$ exists and is equal to
$\rho_A(S)$. Then the definition of $\sigma_A(S)$ implies that the
strict annular density $\sigma_A(S)$ exists and is also equal to
$\rho_A(S)$.  However, as Example~\ref{exmp:even} below shows, it is
possible that $\sigma_A(S)$ exists while $\rho_A(S)$ does not. Thus
there are reasonable situations where the parity of the radius of a
sphere or a ball affects the outcome when measuring the relative size
of a subset of a free group, and annular density turns out to be a
more suitable and relevant quantity.  This is the case when we
consider a subset of $\mathbb Z^k$ and the full preimage of this
subset in $F$ under the abelianization map. Moreover, annular density
and its ``close relatives'' also make sense from the computational
prospective. A typical experiment for generating a ``random'' element
in a ball $B(n)$ of radius $n$ in $F(A)$ might proceed as follows.
First choose a uniformly random integer $m\in [0,n]$ and then choose a
uniformly random element $x$ from the $m$-sphere in $F(A)$ via a
simple non-backtracking random walk of length $m$. It is easy to see
that this experiment, while very natural, does not correspond to the
uniform distribution on $B(n)$. For example, if $F$ has rank $k=2$,
then for the uniform distribution on $B(n)$ the probability that the
element $x$ has length $n$ is approximately $\frac{2}{3}$ for large
$n$ while in our experiment described above this probability is
$1/(n+1)$. In fact if $w\in B(n)$ then the probability of choosing the
element $w$ in the above experiment is $\frac{1}{(n+1)\#S(m)}$ where
$m$ is the freely reduced length of $w$ and where $S(m)$ is the sphere
of radius $m$ in $F$. Thus if $X\subseteq F$ then the probability of
choosing an element of $X$ in the above experiment is
\[
\frac{1}{n+1}\sum_{m=0}^n \frac{\#(X\cap S(m))}{\#S(m)}.
\]
If in our experiment we choose an element of $S(n-1)\cup S(n)$ by
first randomly and uniformly choosing $m\in\{n-1,n\}$ and then
choosing a uniformly random element of $S(m)$, then for a subset $X$
of $F$ the probability of picking an element of $X$ is
\[
\frac{1}{2}\big(\frac{\#(X\cap S(n-1))}{\#S(n-1)}+\frac{\#(X\cap
  S(n))}{\#S(n)}\big),
\]
and the formulas from the definition of annular density appear.

For most of the cases where one can actually compute the asymptotic
density of the set of elements in a free group having some natural
algebraic property, this set turns out to be either generic or
negligible. (Of course, a subset is negligible if and only if its
complement is generic.) The following subsets are known to be
negligible in a free group $F=F(A)$ of rank $k\ge 2$, both in the
sense of asymptotic and annular densities: the set of all proper
powers~\cite{AO}, a finite union of conjugacy classes, a subgroup of
infinite index~\cite{Woess}, a finite union of automorphic orbits
(e.g. the set of all primitive elements)~\cite{Goldstein,KSS}, the set
of all elements whose cyclically reduced forms are not automorphically
minimal~\cite{KSS}, the union of all proper free factors of $F$ (this
follows from results of \cite{St99} and \cite{BV,BMS,KSS}).  Examples
of generic sets, again in the sense of both the asymptotic and the
annular densities, include: the set of all words whose symmetrizations
satisfy the $C'(1/6)$ small cancellation condition~\cite{AO}, the set
of words with nontrivial images in the abelianization of
$F(A)$~\cite{Woess} and the set of elements of $F(A)$ with cyclic
stabilizers in $Aut(F(A))$~\cite{KSS}. It is therefore interesting to
find examples of natural properties of elements of free groups which are
``intermediate'' in the sense that they have density different
from either $0$ or $1$. In this article we show that being a test
element in the free group of rank two is such an example.

\begin{conv}[The abelianization map]
Recall that $F$ is a free group of rank $k\ge 2$ with free basis
$A=\{a_1,\dots, a_k\}$. We identify $\mathbb Z^k$ with the
abelianization of $F$ where the abelianization homomorphism
$\alpha: F\to \mathbb Z^k$ is given by $a_i\mapsto e_i$,
$i=1,\dots, k$. We also denote $\alpha(w)$ by
$\overline{w}$.
\end{conv}

It is easy to construct an example of a subset $H$ of $F$ such that
the annular density of $H$ in $F$ and the asymptotic density of
$\alpha(H)$ in $\mathbb Z^k$ are different. For instance, let
$F=F(a,b)$ and consider the subgroup $H=\langle a,b[a,b]^3\rangle\le
F$. Then $H$ has infinite index in $F$ and hence has both asymptotic
and annular density $0$ in $F$. On the other hand,
$\alpha(H)=\alpha(F)=\mathbb Z^2$ has asymptotic density $1$ in
$\mathbb Z^k$ with respect to any length function on $\mathbb Z^2$.
\begin{exmp}\label{exmp:even}
  Let $F=F(a,b)$, where $A=\{a,b\}$, be free of rank two. Let $\alpha:
  F\to \mathbb Z^2$ be the abelianization map. Note that for any $w\in
  F$ the length $|w|_A$ and $||\alpha(w)||_1$ have the same parity,
  since $||\alpha(w)||_1=|w_a|+|w_b|$, where $w_a, w_b$ are the
  exponent sums on $a$ and $b$ in $w$. Let $S=\{z\in \mathbb Z^2:
  ||z||_1\text{ is even }\}$ and let $\widetilde
  S:=\alpha^{-1}(S)\subseteq F$. Then $\widetilde S=\{w\in F: |w|_A
  \text{ is even }\}$. It is not hard to see that the strict
  asymptotic density of $S$ in $\mathbb Z^2$, with respect to
  $||.||_p$ for any $1\le p\le \infty$, exists and is equal to $1/2$.
  
  Since $\widetilde S$ is exactly the union of all spheres of even
  radii in $F$, and the ratio of the sizes of spheres of radius $n$
  and $n-1$ is equal to $3$, it follows that the limits $\lim_{n\to
    \infty}\frac{\#\{w\in \widetilde S: |w|_A=n\}}{\#\{w\in F:
    |w|=n\}}$ and $\lim_{n\to \infty}\frac{\#\{w\in \widetilde S:
    |w|_A\le n\}}{\#\{w\in F: |w|\le n\}}$ do not exist. However, it
  is easy to see that for every $n\ge 1$
\[
\frac{\#\{w\in \widetilde S:|w|_A=n-1\}}{\#\{w\in F: |w|_A=
  n-1\}}+\frac{\#\{w\in \widetilde S: |w|_A=n\}}{\#\{w\in F:
  |w|_A=n\}}=1,
\]
and therefore $\sigma_A(\widetilde S)=\frac{1}{2}$. Thus although the
strict asymptotic density of $\widetilde S\subseteq F$ does not exist,
the strict annular density does exist and is equal to the strict
asymptotic density of $S\subseteq \mathbb Z^2$. More examples of a
similar nature are discussed in Remark~1.8 of \cite{PR}.
\end{exmp}

Example~\ref{exmp:even} demonstrates why the notion of annular density
is suitable for working with subsets of free groups, while asymptotic
density is more suitable for subsets of free abelian groups.
Geometrically, this difference comes from the fact that free abelian
groups are amenable with balls forming a Folner sequence, while free
groups are non-amenable.

Although the counting occurs in very different places, it is
interesting to ask how the asymptotic density of a subset
$S\subseteq\mathbb Z^k$, with respect to some natural length function,
and the annular density of its full preimage $\alpha^{-1}(S)$ in $F$
are related.  We shall see that there is a reasonable assumption about
the set $S$ which guarantees that the two densities are actually
equal.

To do this we need to understand the image of the uniform distribution
on the sphere of radius $n$ in $F$ under the abelianization map
$\alpha$.  There is an explicit formula for the size of the preimage
of an element, and there is also a Central Limit Theorem saying that,
when appropriately normalized, the distribution converges to a normal
distribution. The methods of~\cite{Riv99} also give a \emph{Local
  Limit Theorem} showing that, when working with width-two spherical
shells in a free group, the densities of the image distributions in
$\mathbb Z^k$ converge to a normal density. Such a result was later
also shown (by rather different methods, and in greater generality) by
Richard Sharp in~\cite{Sharp}. Recently Petridis and Risager~\cite{PR}
obtained a similar Local Limit Theorem for counting conjugacy classes
rather than elements of $F$. On the face of it, studying the annular
density of the set $\alpha^{-1}(S)$ in $F$ presents new challenges.
The central limit theorem by itself seems too crude a tool and \emph{a
  priori} it would appear that one would need very sharp error bounds
in the local limit theorem.  Nevertheless, we produce a short argument
solving this problem where one of the key ingredients is the
ergodicity of the $SL(k,\mathbb Z)$-action on $\mathbb R^k$.

We can now state our main result:

\begin{theor}\label{A'}
  Let $F=F(A)$ be a free group of rank $k\ge 2$ with free basis
  $A=\{a_1,\dots, a_k\}$ and let $\alpha:F\to \mathbb Z^k$ be the
  abelianization homomorphism.

  Let $S\subseteq \mathbb Z^k$ be an $SL(k, \mathbb Z)$-invariant
  subset and put $\widetilde{S} =\alpha^{-1}(S)\subseteq F$.
  
  Then
\begin{enumerate}
  
\item For every $1\le p\le \infty$ the strict asymptotic density
  $\rho_p(S)$ exists and, moreover, for every $1\le p\le \infty$ we
  have $\rho_p(S)=\rho_\infty(S)$.
  
\item The strict annular density $\sigma_{A}(\widetilde{S})$ exist
  and, moreover, $\sigma_{A}(\widetilde{S})=\rho_{\infty}(S)$.
  
  That is,

\begin{gather*}
  \lim_{n\to\infty} \frac{1}{2}\big(\frac{\gamma_A(n-1,\{w\in F:
    \alpha(w)\in S \})}{\gamma_A(n-1,F)}+\frac{\gamma_A(n,\{w\in F:
    \alpha(w)\in S \}}{\gamma_A(n,F)}\big)=\\
  \lim_{n\to\infty} \frac{\#\{z: z\in \mathbb Z^k, ||z||_\infty\le n,
    \text{ and } z\in S \}}{\#\{z: z\in \mathbb Z^k, ||z||_\infty \le
    n\}}.
\end{gather*}
\end{enumerate}
\end{theor}

The requirement that $S$ be $SL(k,\mathbb Z)$-invariant essentially
says that the subset $S$ of $\mathbb Z^k$ is defined in ``abstract"
group-theoretic terms, not involving the specific choice of a free
basis for $\mathbb Z^k$. Note that Proposition~\ref{z^k} below gives
an explicit formula for $\rho_{\infty}(S)$ in Theorem~\ref{A'}.

Our main application of Theorem~\ref{A'} concerns the case where $S$
is the set of all ``visible" points in $\mathbb Z^k$. A nonzero point
$z$ of $\mathbb Z^k$ is called \emph{visible} if the greatest common
divisor of the coordinates of $z$ is equal to $1$.  This terminology
is standard in number theory~\cite{Ru} and reflects the fact that if
$z$ is visible then the line segment between the origin and $z$ does
not contain any other integer lattice points.  For a nonzero point
$z\in \mathbb Z^k$ being visible is also equivalent to $z$ not being a
proper power in $\mathbb{Z}^k$, that is, to $z$ generating a maximal
cyclic subgroup of $\mathbb{Z}^k$.  More generally, if $t\ge 1$ is an
integer, we will say that $z\in \mathbb Z^k$ is \emph{$t$-visible} if
$z=z_1^t$ for some visible $z_1\in \mathbb Z^k$, that is, if the
greatest common divisor of the coordinates of $z$ is equal to $t$.

We want to ``lift'' this terminology to free groups.

\begin{defn}[Visible elements in free groups]
  Let $F=F(A)$ be a free group of rank $k\ge 2$ with free basis
  $A=\{a_1,\dots, a_k\}$ and let $\alpha:F\to \mathbb Z^k$ be the
  abelianization homomorphism, that is, $\alpha(a_i)=e_i\in \mathbb
  Z^k$.  We say that an element $w\in F$ is \emph{visible} if
  $\alpha(w)$ is visible in $\mathbb Z^k$. Let $V$ be the set of
  visible elements of $F$. Similarly, for an integer $t\ge 1$ an
  element $w\in F$ is \emph{$t$-visible} if $\alpha(w)$ is $t$-visible
  in $\mathbb Z^k$. We use $V_t$ to denote the set of all $t$-visible
  elements of $F$ and we use $U_t$ to denote the set of all
  $t$-visible elements of $\mathbb Z^k$.
\end{defn}
Note that $V=V_1$ and that for every $t\ge 1$ the definition of $V_t$
does not depend on the choice of the free basis $A$ of $F$.

The following proposition giving the asymptotic density of the set of
$t$-visible points in $\mathbb Z^k$ in terms of the Riemann
zeta-function is well-known in number theory~\cite{Chr}.

\begin{prop}\label{u}
  For any integer $t\ge 1$ we have
\[
\rho_\infty(U_t)=\frac{1}{t^k\zeta(k)}.
\]
\end{prop}

The case $k =2$ and $t = 1$ of Proposition~\ref{u} was proved by
Mertens in 1874~\cite{Mer}. (See also Theorem~332 of the classic book
of Hardy and Wright~\cite{HW}.) Recall that
$\zeta(k)=\sum_{n=1}^{\infty} \frac{1}{n^k}$ and, in particular,
$\zeta(2)=\frac{\pi^2}{6}$.

It is therefore natural to investigate the asymptotic density of the
set of visible elements in $F$.  As a direct corollary of
Theorem~\ref{A'}, of Proposition~\ref{prop:compare} below and of
Proposition~\ref{u} we obtain:

\begin{theor}\label{A}
  Let $F=F(A)$ be a free group of rank $k\ge 2$ with free basis
  $A=\{a_1,a_2,\dots, a_k\}$. Let $t\ge 1$ be an integer. Then the
  strict annular density $\sigma_A(V_t)$ exists and
  \[\sigma_A(V_t)=\frac{1}{t^k\zeta(k)}.\]
  
  Moreover, in this case
\begin{gather*}
  0<\frac{4k-4}{(2k-1)^2t^k\zeta(k)}\le \liminf_{n\to\infty}
  \frac{\rho_A(n,V_t)}{\rho_A(n,F)}\le\\ \limsup_{n\to\infty}
  \frac{\rho_A(n,V_t)}{\rho_A(n,F)}\le
  1-\frac{4k-4}{(2k-1)^2}\big(1-\frac{1}{t^k\zeta(k)}\big)<1.
\end{gather*}

\end{theor}
A result similar to Theorem~\ref{A} for counting conjugacy classes in
$F_k$ with primitive images in $\mathbb Z^k$ has been recently
independently obtained by Petridis and Risager~\cite{PR}.

For the case of the free group of rank two we compute the two
``spherical densities'' for the density of the set $V_1$,
corresponding to even and odd $n$ tending to infinity:

\begin{theor}\label{spher}
  Let $k=2$. We have
\[
\lim_{m\to\infty} \frac{\gamma_A(2m,
  V_1)}{\gamma_A(2m,F)}=\frac{2}{3\zeta(2)}=\frac{4}{\pi^2}
\]
and
\[
\lim_{m\to\infty} \frac{\gamma_A(2m-1,
  V_1)}{\gamma_A(2m-1,F)}=\frac{8}{\pi^2}.
\]
\end{theor}
Theorem~\ref{A} implies that $\sigma_A(V_1)=\frac{6}{\pi^2}$ in this
case. Since the two limits in Theorem~\ref{spher} are different, the
statements of Theorem~\ref{A'} and Theorem~\ref{A} cannot be
substantially improved. This fact underscores the conclusion that
annular density is the right kind of notion for measuring the sizes of
subsets of free groups, where the abelianization map is concerned.

We apply Theorem~\ref{A} to compute the annular density of test
elements in a free group of rank two.  Recall that an element $g\in G$
is called a \emph{test element} if every endomorphism of $G$ fixing
$g$ is actually an automorphism of $G$.  It is easy to see that for
two conjugate elements $g_1,g_2\in G$ the element $g_1$ is a test
element if and only if $g_2$ is also a test element and thus the
property of being a test element depends only on the conjugacy class
of an element $g$.  The notion of a test element was introduced by
Shpilrain \cite{Sh94} and has since become a subject of active
research both in group theory and in the context of other algebraic
structures such as polynomial algebras and Lie algebras. (See, for
example, \cite{FRSS,Iv,Lee,MSY,MUY,OT,Ro,Sh95,Ti}.) It turns out that
studying test elements in a particular group $G$ produces interesting
information about the automorphism group of $G$.

Here we prove:

\begin{theor}\label{B}
  Let $F=F(a,b)$ be a free group of rank two with free basis
  $A=\{a,b\}$.  Then for the set $\mathcal T$ of all test elements in
  $F$ the strict annular density exists and
\[ \sigma_A(\mathcal T)=1-\frac{6}{\pi^2}.\] Moreover,
\begin{gather*}
0<\frac{4}{9}(1-\frac{6}{\pi^2})\le \liminf_{n\to\infty}
\frac{\rho_A(n,\mathcal T)}{\rho_A(n,F)}\le\\ \limsup_{n\to\infty}
\frac{\rho_A(n,\mathcal T)}{\rho_A(n,F)}\le 1-\frac{8}{3\pi^2}<1.
\end{gather*}

\end{theor}

By a result of Turner \cite{Tu}, an element of $F$ is a test elements
if and only if it does not belong to a proper retract of $F$.
Therefore Theorem~\ref{B} implies that the strict annular density of
the union of all proper retracts of $F(a,b)$ is $\frac{6}{\pi^2}$. In
Theorem~\ref{B} above we have $\frac{4}{9}(1-\frac{6}{\pi^2})\approx
0.1742$, $1-\frac{6}{\pi^2}\approx .3920$ and
$1-\frac{8}{3\pi^2}\approx 0.7298$. Thus Theorem~\ref{B} shows that
being a test element is an ``intermediate'' property in the free group
of rank two. More generally, Theorem~\ref{A} implies that, for every
$k \ge 2$ and for $N$ sufficiently large, the set $V_N$ of $N$-visible
elements in $F_k$ has strictly positive annular density arbitrarily
close to $0$ while the set $S_N = \bigcup_{t = 1}^{N} V_t$ has annular
density less than but arbitrarily close to $1$. It is well-known that
every positive rational number has an ``Egyptian fraction''
representation as a finite sum of distinct terms of the form
$\frac{1}{n}$. It does not seem clear what values can be obtained as
finite sums of distinct terms of the form $\frac{1}{n^2}$ and there
are excluded intervals. In particular, if such a sum uses $1$, it is
at least $1$ while if $1$ is not used then the sum is at most
$\frac{\pi ^2}{6} -1$. Multiplying by the scale factor $\frac{1}{\zeta
  (2)}$, we see that we cannot obtain a annular density in the open
interval $(1 - \frac{6}{\pi^2}, \frac{6}{\pi^2} )$ by taking a finite
union of the sets $V_t$ (Proposition~\ref{z^k} below shows that the
same is true for infinite unions). It is interesting to note that the
probabilities of being a test element or of not being a test element
are the boundary points of this excluded interval.

Nathan Dunfield and Dylan Thurston recently proved~\cite{DT} that for
a two-generator one-relator group being free-by-cyclic is an
intermediate property. While they do not provide an exact value for
the asymptotic density (nor do they prove that either the strict
asymptotic density or the strict annular density exist), they show
that it is strictly between $0$ and $1$. Computer experiments by
Kapovich and Schupp, by Mark Sapir and by Dunfield and Thurston
indicate that in the two-generator case this asymptotic density is
greater than $0.9$.

The authors are grateful to Laurent Bartholdi, John D'Angelo,Iwan
Duursma, Kevin Ford, Steve Lalley, Alexander Ol'shanskii, Yuval Peres,
Yannis Petridis, Alexandru Zaharescu and Andrzej Zuk for very helpful
conversations.

\section{Comparing densities in $\mathbb Z^k$ and in $F_k$}

\begin{conv}
  Throughout this section let $S\subseteq \mathbb Z^k$ be as in
  Theorem~\ref{A'} and let $\delta:=\overline{\rho}_\infty(S)$.
\end{conv}

We can now prove that for every $SL(k,\mathbb Z)$-invariant subset $S$
of $\mathbb Z^k$ the strict asymptotic density ${\rho}_\infty(S)$
exists.  Recall that Proposition~\ref{u} stated in the introduction
gives the precise value of the strict asymptotic density of the set of
all $t$-visible elements in $\mathbb Z^k$. The crucial points of the
proof below are that the complement of an $SL(k,\mathbb Z)$-invariant
is also $SL(k,\mathbb Z)$-invariant and that $\sum_{t=1}^\infty
\frac{1}{t^k\zeta(k)}=1$.

\begin{prop}\label{z^k}
  Let $Y\subseteq Z^k$ be a nonempty $SL(k,\mathbb Z)$-invariant
  subset that does not contain $\mathbf 0\in \mathbb Z^k$. Let $I$ be
  the set of all integers $t\ge 1$ such that there exists a
  $t$-visible element in $Y$. Then
\begin{enumerate}
\item $Y=\cup_{t\in I} U_t$.
  
\item The strict asymptotic density $\rho_\infty(Y)$ exists and
\[
\rho_\infty(Y)=\sum_{t\in I}\rho(U_t)=\sum_{t\in
  I}\frac{1}{t^k\zeta(k)}.
\]
\end{enumerate}
\end{prop}
\begin{proof}
  
  Observe first that
\[
\sum_{t=1}^\infty
\frac{1}{t^k\zeta(k)}=\frac{1}{\zeta(k)}\sum_{t=1}^\infty
\frac{1}{t^k}=\frac{\zeta(k)}{\zeta(k)}=1.
\]

Since $k\ge 2$, two nonzero elements $z,z'\in \mathbb Z^k$ lie in the
same $SL(k,\mathbb Z)$-orbit if and only if the greatest common
divisors of the coordinates of $z$ and of $z'$ are equal.  Thus every
$SL(k,\mathbb Z)$-orbit of a nonzero element of $\mathbb Z^k$ has the
form $U_t$ for some $t\ge 1$. This implies part (1) of
Proposition~\ref{z^k}.

Let $I':=\{t\in \mathbb Z: t\ge 1, t\not\in I\}$. If either $I$ or
$I'$ is finite, part (2) of Proposition~\ref{z^k} follows directly
from proposition~\ref{u}. Suppose now that both $I$ and $I'$ are
infinite and let $Y':=\mathbb Z^k-(Y\cup\{\mathbf 0\})=\cup_{t\in I'}
U_t$.

For every finite subset $J\subseteq I$ let $Y_J:=\cup_{t\in J} U_t$.
Since $Y_J\subseteq Y$, it follows that
\[
\liminf_{n\to\infty}\frac{\#\{z\in Y: ||z||_\infty\le n\}}{\#\{z\in
  \mathbb Z^k: ||z||_\infty\le n\}}\ge \rho_\infty(Y_J)=\sum_{t\in J}
\frac{1}{t^k\zeta(k)}.
  \]
  Since this is true for every finite subset of $I$, we conclude that
\[
\liminf_{n\to\infty}\frac{\#\{z\in Y: ||z||_\infty\le n\}}{\#\{z\in
  \mathbb Z^k: ||z||_\infty\le n\}}\ge \sum_{t\in I}
\frac{1}{t^k\zeta(k)}.
  \]
  
  The same argument applies to the $SL(k,\mathbb Z)$-invariant set
  $Y'$ and therefore:
\[
\liminf_{n\to\infty}\frac{\#\{z\in Y': ||z||_\infty\le n\}}{\#\{z\in
  \mathbb Z^k: ||z||_\infty\le n\}}\ge \sum_{t\in I'}
\frac{1}{t^k\zeta(k)}.
\]
This implies
\begin{gather*}
  1-\liminf_{n\to\infty}\frac{\#\{z\in Y': ||z||_\infty\le
    n\}}{\#\{z\in \mathbb Z^k: ||z||_\infty\le n\}}\le 1-\sum_{t\in
    I'}
  \frac{1}{t^k\zeta(k)}\Rightarrow\\
  \limsup_{n\to\infty}\left(1-\frac{\#\{z\in Y': ||z||_\infty\le
      n\}}{\#\{z\in \mathbb Z^k: ||z||_\infty\le n\}}\right)\le
  1-\sum_{t\in I'}
  \frac{1}{t^k\zeta(k)}\Rightarrow\\
  \limsup_{n\to\infty}\frac{\#\{z\in Y: ||z||_\infty\le n\}}{\#\{z\in
    \mathbb Z^k: ||z||_\infty\le n\}}\le 1-\sum_{t\in I'}
  \frac{1}{t^k\zeta(k)}=\sum_{t\in I} \frac{1}{t^k\zeta(k)}.
\end{gather*}

Hence \[ \lim_{n\to\infty}\frac{\#\{z\in Y: ||z||_\infty\le
  n\}}{\#\{z\in \mathbb Z^k: ||z||_\infty\le n\}}=\sum_{t\in I}
\frac{1}{t^k\zeta(k)},
\]
as required.
\end{proof}

Recall that $S\subseteq \mathbb Z^k$ is an $SL(k,\mathbb Z)$-invariant
subset and that $\delta=\overline{\rho}_\infty(S)$.
Proposition~\ref{z^k} implies that in fact $\delta=\rho_\infty(S)$.

It is well known that if $\Omega\subseteq \mathbb R^k$ is a ``nice''
bounded open set then the Lebesgue measure $\lambda(\Omega)$ can be
computed as

\[
\lambda(\Omega)=\lim_{r\to\infty}\frac{\# (\mathbb Z^k\cap r\Omega
  )}{r^k}.
\]
Here we say that a bounded open subset of $\mathbb R^k$ is ``nice'' if
its boundary is piecewise smooth.

We need a similar formula for counting the points of $S$.  For a real
number $r\ge 1$ and a nice bounded open set $\Omega\subseteq \mathbb
Z^k$ let
\[
\mu_{r,S}(\Omega):=\frac{\#(S\cap r\Omega)}{r^k}.
\]

\begin{prop}\label{prop:zeta}
  For any nice bounded open set $\Omega\subseteq \mathbb Z^k$ we have
\[
\lim_{r\to\infty} \mu_{r,S}(\Omega)=\delta\lambda(\Omega).
\]
\end{prop}

\begin{proof}
  Each $\mu_{r,S}$ can be regarded as a measure on $\mathbb R^k$.  We
  prove the theorem by showing that the $\mu_{r,S}$ weakly converge to
  $\delta\lambda$ as $r\to\infty$, where $\lambda$ is the Lebesgue
  measure.

  By Helly's theorem there exists a sequence $(r_i)_{i=1}^{\infty}$
  with $\displaystyle\lim_{i\to\infty} r_i=\infty$ such that the
  sequence $\mu_{r_1,S},\mu_{r_2,S}, \dots$ is weakly convergent to
  some limiting measure.  We now show that for every such convergent
  subsequence of $\mu_{r_i,S}$ the limiting measure is indeed equal to
  $\delta\lambda$, where $\lambda$ is the Lebesgue measure.

  Indeed, suppose that $\sigma=(r_i)_{i=1}^{\infty}$ is a sequence
  with $\displaystyle\lim_{i\to\infty} r_i=\infty$ such that the
  sequence $\mu_{r_i,S}$ converges to the limiting measure
  $\displaystyle\mu_\sigma=\lim_{i\to\infty} \mu_{r_i,S}$.  Every
  $\mu_{r_i,S}$ is invariant with respect to the natural $SL(k,\mathbb
  Z)$-action since this action preserves the set $S$ and also commutes
  with homotheties of $\mathbb R^k$ centered at the origin.  Therefore
  the limiting measure $\mu_\sigma$ is also $SL(k,\mathbb
  Z)$-invariant.

  Moreover, the measures $\mu_{r,S}$ are dominated by the measures
  $\lambda_r$ defined as $\lambda_r(\Omega)=\frac{\# (\mathbb Z^k\cap
    r\Omega )}{r^k}$. Since, as observed earlier, the measures
  $\lambda_r$ converge to the Lebesgue measure $\lambda$, it follows
  that $\mu_\sigma$ is absolutely continuous with respect to
  $\lambda$.  It is known that the natural action of $SL(k,\mathbb Z)$
  on $\mathbb R^k$ is ergodic with respect to $\lambda$. (See, for
  example, Zimmer's classic monograph~\cite{Zim}.)  Therefore
  $\mu_\sigma$ is a constant multiple $c\lambda$ of $\lambda$. The
  constant $c$ can be computed explicitly for a set such as the open
  unit ball $B$ in the $||.||_\infty$ norm on $\mathbb R^k$ defining
  the length function $\ell$ on $\mathbb Z^k$.

  By assumption we know that
  \[
  \lim_{r\to\infty} \frac{\#\{z\in \mathbb Z^k: z\in S\cap
    rB\}}{\#\{z\in \mathbb Z^k: z\in rB\}}=\delta.
  \]
  
  We also have
\[
\lim_{r\to\infty} \frac{\#\{z\in \mathbb Z^k: z\in
  rB\}}{r^k}=\lambda(B)
\]
and hence
\[
\lim_{r\to\infty} \frac{\#\{z\in \mathbb Z^k: z\in S\cap
  rB\}}{r^k}=\delta\lambda(B).
\]
Therefore $c=\delta$ and $\mu_\sigma=\delta\lambda$. The above
argument in fact shows that every convergent subsequence, with
$r\to\infty$, of $\mu_{r,S}$ converges to $\delta\lambda$ and
therefore $\displaystyle\lim_{r\to\infty} \mu_{r,S} =\delta\lambda$.

\end{proof}
\begin{rem}\label{rem:p}
  Let $1\le p\le \infty$. Then the open unit ball in $\mathbb R^k$
  with respect to $||.||_p$ is "nice".  Proposition~\ref{prop:zeta},
  applied to $\Omega$ being this ball, implies that
  $\rho_{p}(S)=\rho_\infty(S)=\delta$.
\end{rem}

\begin{conv} As always, $F=F(a_1,\dots, a_k)$ is the free group of rank
  $k\ge 2$ with free basis $A=\{a_1,\dots, a_k\}$ and $\alpha:F\to
  \mathbb Z^k$ is the abelianization homomorphism sending $a_i$ to
  $e_i$ in $\mathbb Z^k$. We will denote $\alpha(w)$ by $\overline w$.
  For $n\ge 1$, $B_F(n)$ denotes the set of all $w\in F$ with
  $|w|_A\le n$. Also, for a point $x=(x_1,\dots,x_k)\in \mathbb R^k$
  we denote by $||x||$ the $||.||_2$-norm of $x$, that is
  $||x||=\sqrt{\sum_{i=1}^k x_i^2}$.
\end{conv}

\begin{notation}\label{not:pn}
  For an integer $n\ge 1$ and a point $x\in \mathbb R^k$ let
  \[
p_n(x)=\frac{\gamma_A(n-1,\{f\in F:
\alpha(f)=x\sqrt{n}\})}{2\gamma_A(n-1,F)}+ \frac{\gamma_A(n,\{f\in
F: \alpha(f)=x\sqrt{n}\})}{2\gamma_A(n,F)}.
\]
\end{notation}
Thus $p_n$ is a distribution supported on finitely many points of
$\frac{1}{\sqrt{n}} \mathbb Z^k$.

We need the following facts about the sequence of distributions $p_n$.
Of these the most significant is part (2) which is a local limit
theorem in our context. It was obtained by Rivin~\cite{Riv99} and,
independently and via different methods, by Sharp~\cite{Sharp}
(specifically, we use Theorem~1 of \cite{Sharp} for part (2) of
Proposition~\ref{prop:limit} below).

\begin{prop}\label{prop:limit}\cite{Riv99,Sharp}
  Let $k\ge 2$ and let $p_n$ be as above. Then:

\begin{enumerate}
\item The sequence of distributions $p_n$ converges weakly to a normal
  distribution $\mathfrak{N}$, with density $\mathfrak{n}$.
  
\item We have
\[
\sup_{x\in \mathbb Z^k/\sqrt{n}} |p_n(x)n^{k/2}-\mathfrak{n}(x)|\to 0
\text{ as } n\to \infty.
\]

\item We have
\[
\lim_{c\to\infty} \sum \{ p_n(x): x\in \mathbb Z^k/\sqrt{n} \text{ and
} ||x||\ge c\}=0.
\]

\end{enumerate}
\end{prop}

\begin{thm}\label{thm:vis-free}
  Let $\Omega\subseteq \mathbb R^k$ be a nice bounded open set. Then
\[
\lim_{n\to\infty} \sum_{x\in S\cap \sqrt{n}\Omega}
p_n(x/\sqrt{n})=\delta\mathfrak{N}(\Omega).
\]
\end{thm}
\begin{proof}
  We have
\begin{gather*}
  \sum_{x\in \mathbb Z_t^k\cap \sqrt{n}\Omega}
  p_n(x/\sqrt{n})=\sum_{y\in \frac{1}{\sqrt{n}} S\cap \Omega}
  p_n(y)=\\
  {n}^{-k/2} \sum_{y\in \frac{1}{\sqrt{n}}S\cap
    \Omega}\mathfrak{n}(y)\ \ +\ \ {n}^{-k/2} \sum_{y\in
    \frac{1}{\sqrt{n}}S\cap \Omega} ({n}^{k/2}
  p_n(y)-\mathfrak{n}(y)).
\end{gather*}

The local limit theorem in part (2) of Proposition~\ref{prop:limit}
tells us that, as $n\to\infty$, each summand ${n}^{k/2}
p_n(y)-\mathfrak{n}(y)$ of the second sum in the last line of equation
above converges to zero and hence so does their Cesaro mean.
Proposition~\ref{prop:zeta} implies that, as $n\to\infty$, the first
summand $\displaystyle {n}^{-k/2} \sum_{y\in \frac{1}{\sqrt{n}}S\cap
  \Omega}\mathfrak{n}(y)$ converges to
 \[
 \delta \int_\Omega \mathfrak{n} d\lambda =
 \delta\mathfrak{N}(\Omega).
\]
\end{proof}

We can now compute the strict asymptotic density of $\widetilde{S}
=\alpha^{-1}(S)$ in $F$ and obtain Theorem~\ref{A'}.

\begin{proof}[Proof of Theorem~\ref{A'}]
  
  Recall that $S\subseteq \mathbb Z^k$ is an $SL(k,\mathbb
  Z)$-invariant set and that $\delta=\overline\rho_\infty(S)$.
  Proposition~\ref{z^k} implies that in fact $\rho_\infty(S)$ exists
  and $\delta=\rho_\infty(S)$. Moreover, as we have seen in
  Remark~\ref{rem:p}, for every $1\le p\le \infty$ the strict
  asymptotic density $\rho_p(S)$ exists and
  $\rho_p(S)=\delta=\rho_\infty(S)$. This proves part (1) of
  Theorem~\ref{A'}.
  
  To prove part (2) of Theorem~\ref{A'} we need to establish that the
  strict annular density $\sigma_A(S)$ exists and that
  $\sigma_A(S)=\delta$.
  
  For $c>0$ denote $\Omega_c:=\{x\in \mathbb R^k: ||x||<c\}$. Then
  $\displaystyle\lim_{c\to\infty} \mathfrak{N}(\Omega_c)=1$.  Let
  $\epsilon>0$ be arbitrary. Choose $c>0$ such that
  \[|\mathfrak{N}(\Omega_c)-1|\le \epsilon/3\] and such that
\[
\lim_{n\to\infty}\sum \{ p_n(x): x\in \mathbb Z^k/\sqrt{n} \text{ and
} ||x||\ge c\}\le \epsilon/6.
\]

By Theorem~\ref{thm:vis-free} and the above formula there is some
$n_0\ge 1$ such that for all $n\ge n_0$ we have
\[
\left| \sum_{x\in S\cap \sqrt{n}\Omega_c}
  p_n(x/\sqrt{n})-\delta\mathfrak{N}(\Omega_c) \right| \le \epsilon/3
\]
and
\[
\sum \{ p_n(x): x\in \mathbb Z^k/\sqrt{n} \text{ and } ||x||\ge c\}\le
\epsilon/3.
\]

Let \[Q(n):=\frac{\gamma_A(n-1, \{w\in F: \overline{w}\in
  S\})}{2\gamma_A(n-1,F)}+\frac{\gamma_A(n, \{w\in F: \overline{w}\in
  S\})}{2\gamma_A(n,F)}.\]

For $n\ge n_0$ we have
\begin{gather*}
  Q(n)=\\
   \frac{\#\{w\in F: \overline{w}\in
    S,  |w|_A=n-1 \text{ and } ||\overline
    w||<c\sqrt{n}\}}{2\gamma_A(n-1,F)} +\\ \frac{\#\{w\in F: \overline{w}\in
    S,  |w|_A=n \text{ and } ||\overline
    w||<c\sqrt{n}\}}{2\gamma_A(n,F)}+\\
\frac{\#\{w\in F: \overline{w}\in
    S,  |w|_A=n-1 \text{ and } ||\overline
    w||\ge c\sqrt{n}\}}{2\gamma_A(n-1,F)} +\\ \frac{\#\{w\in F: \overline{w}\in
    S,  |w|_A=n \text{ and } ||\overline
    w||\ge c\sqrt{n}\}}{2\gamma_A(n,F)}=\\
  \sum_{x\in S\cap \sqrt{n}\Omega_c} p_n(x/\sqrt{n})\ \ +\ \sum_{x\in
    S \cap (\mathbb R^k-\sqrt{n}\Omega_c)} p_n(x/\sqrt{n})
\end{gather*}

In the last line of the above equation, the first sum differs from
$\delta\mathfrak{N}(\Omega_c)$ by at most $\epsilon/3$ since $n\ge
n_0$ and the second sum is $\le \epsilon/3$ by the choice of $c$ and
$n_0$. Therefore, again by the choice of $c$, we have
$|Q(n)-\delta|\le \epsilon$. Since $\epsilon>0$ was arbitrary, this
implies that $\displaystyle\lim_{n\to\infty} Q(n)=\delta$, as claimed.
\end{proof}

The following observation shows how to estimate the asymptotic density
in terms of the annular density.
\begin{prop}\label{prop:compare}
  Let $Y\subseteq F$ be a subset such that the strict annular density
  $\delta=\sigma_A(Y)$ exists. Then
\[
\frac{4k-4}{(2k-1)^2}\delta\le
\liminf_{n\to\infty}\frac{\rho_A(n,S)}{\rho_A(n,F)}\le
\limsup_{n\to\infty}\frac{\rho_A(n,S)}{\rho_A(n,F)}\le
1-\frac{4k-4}{(2k-1)^2}(1-\delta).
\]
In particular, if $0<\delta<1$ then
\[
0<\liminf_{n\to\infty}\frac{\rho_A(n,S)}{\rho_A(n,F)}\le
\limsup_{n\to\infty}\frac{\rho_A(n,S)}{\rho_A(n,F)}<1.
\]
\end{prop}
\begin{proof}
  Note that for $n\ge 1$ we have $\gamma_A(n,F)=2k(2k-1)^n$ and that,
  up to an additive constant, $\rho_A(n,F)=\frac{k}{k-1}(2k-1)^n$.
  Denote $a_n=\gamma_A(n,Y)$.  We have
\[
\delta=\lim_{n\to\infty}\frac{1}{2}\big(
\frac{a_{n-1}}{2k(2k-1)^{n-2}}+\frac{a_{n}}{2k(2k-1)^{n-1}}\big)=
\frac{1}{2}\lim_{n\to\infty}\frac{a_{n-1}\frac{2k-1}{2k-2}+a_n\frac{1}{2k-2}}{\frac{k}{k-1}(2k-1)^{n-1}}
\]
Therefore
\begin{gather*}
  \liminf_{n\to\infty}\frac{\rho_A(n,Y)}{\rho_A(n,F)}=\liminf_{n\to\infty}\frac{a_1+\dots
    +a_n}{\frac{k}{k-1}(2k-1)^n}\ge\liminf_{n\to\infty}\frac{a_{n-1}+a_n}{\frac{k}{k-1}(2k-1)^n}=\\
  \frac{2k-2}{2k-1}\liminf_{n\to\infty}\frac{a_{n-1}\frac{2k-1}{2k-2}+a_n\frac{2k-1}{2k-2}}{\frac{k}{k-1}(2k-1)^n}\ge
  \frac{2k-2}{2k-1}\liminf_{n\to\infty}\frac{a_{n-1}\frac{2k-1}{2k-2}+a_n\frac{1}{2k-2}}{\frac{k}{k-1}(2k-1)^n}=\\
  \frac{4k-4}{(2k-1)^2}\liminf_{n\to\infty}\frac{1}{2}\cdot
  \frac{a_{n-1}\frac{2k-1}{2k-2}+a_n\frac{1}{2k-2}}{\frac{k}{k-1}(2k-1)^{n-1}}=\frac{4k-4}{(2k-1)^2}\delta.
\end{gather*}

Applying the same argument to the set $F-Y$, we get
\[
\liminf_{n\to\infty}\frac{\rho_A(n,F-Y)}{\rho_A(n,F)}\ge
\frac{4k-4}{(2k-1)^2}(1-\delta).
\]
Therefore
\[
\limsup_{n\to\infty}\frac{\rho_A(n,Y)}{\rho_A(n,F)}=1-\liminf_{n\to\infty}\frac{\rho_A(n,F-Y)}{\rho_A(n,F)}\le
1-\frac{4k-4}{(2k-1)^2}(1-\delta).
\]
\end{proof}

\section{Spherical densities}

In this section we will prove Theorem~\ref{spher} from the
Introduction and, for the case of $k=2$, compute the ``spherical
densities''
\[
\lim_{m\to\infty} \frac{\gamma_A(2m, V_1)}{\gamma_A(2m,F)}
\]
and
\[
\lim_{m\to\infty} \frac{\gamma_A(2m-1, V_1)}{\gamma_A(2m-1,F)}.
\]
This is done by computing the ``spherical densities'' for the set
$V_1(e)\subseteq F$ consisting of all points of $V_1$ of even length
and comparing it with the strict asymptotic density of the set
$U_1(e)\subseteq \mathbb Z^k$ of all elements of $\mathbb Z^k$ of even
$||.||_\infty$-length. The key point is that for $n=2m$ we have
$\gamma_A(2m-1, V_1(e))=0$ and $\gamma_A(2m,V_1)=\gamma_A(2m,
V_1(e))$. Therefore for the quantities from the definition of annular
density we have
\[
\frac{1}{2}\big(\frac{\gamma_A(2m,V_1(e))}{\gamma_A(2m,F)}+\frac{\gamma_A(2m-1,
  V_1(e))}{\gamma_A(2m-1,F)}\big)=\frac{\gamma_A(2m,
  V_1(e))}{2\gamma_A(2m,F)}=\frac{\gamma_A(2m, V_1)}{2\gamma_A(2m,F)}.
\]
This allows us to essentially repeat the proof of Theorem~\ref{A'},
applied to the sets $U_1(e)$ and $V_1(e)=\alpha^{-1}(U_1(e))$, except
that instead of ergodicity of the action of $SL(k,\mathbb Z)$ we use
the ergodicity of the action on $\mathbb R^k$ of a congruence subgroup
of $SL(k,\mathbb Z)$ that leaves $U_1(e)$ invariant.

\begin{conv}
  We say that an element $z=(z_1,\dots,z_k)\in \mathbb Z^k$ is
  \emph{even} if $||z||_1=|z_1|+\dots+|z_k|$ is even and that $z$ is
  \emph{odd} if $||z||_1$ is odd. Similarly, $w\in F$ is \emph{even}
  if $|w|_A$ is even and $w\in F$ is \emph{odd} if $|w|_A$ is odd.
  Note that $w\in F$ is even if and only if $\alpha(w)\in \mathbb Z^k$
  is even.
  
  Let $G_k$ be the set of all $M\in SL(k,\mathbb Z)$ such that $M=I_k$
  in $SL(k,\mathbb Z/2\mathbb Z)$. Thus $G_k$ is a finite index
  subgroup of $SL(k,\mathbb Z)$ also known as the $2$-congruence
  subgroup. Denote by $\mathbb Z^k(e)$ the set of all even elements in
  $\mathbb Z^k$. Also denote $U_1(e):=U_1\cap \mathbb Z^k(e)$. Observe
  that $\mathbb Z^k(e)$ and $U_1(e)$ are $G_k$-invariant and that
  $V_1(e)=\alpha^{-1}(U_1(e))$. (The actual set-wise stabilizer of
  $\mathbb Z^k(e)$ in $SL(k,\mathbb Z)$ contains $G_k$ as a subgroup
  of finite index.)
\end{conv}

\begin{prop}\label{prop:zeta1}
Let $S\subseteq \mathbb Z^k$ be a subset such that
$\delta=\rho_\infty(S)$ exists and such that $S$ is
$G_k$-invariant. Then for every bounded nice open set
$\Omega\subseteq \mathbb R^k$ we have

\[
\lim_{r\to\infty} \mu_{r,S}(\Omega)=\delta\lambda(\Omega).
\]
\end{prop}
\begin{proof}
  The proof is the same as for Proposition~\ref{prop:zeta}. The only
  difference is that instead of ergodicity of the $SL(k,\mathbb
  Z)$-action on $\mathbb R^k$ we use ergodicity of the $G_k$-action on
  $\mathbb R^k$ with respect to the Lebesgue measure (see \cite{Zim}
  for the proof of this ergodicity).
\end{proof}

Let $p_n(x)$ be defined exactly as in Notation~\ref{not:pn}.

\begin{thm}\label{thm:vis-free1}
  Let $\Omega\subseteq \mathbb R^k$ be a nice bounded open set.  Let
  $S\subseteq \mathbb Z^k$ be a $G_k$-invariant subset such that
  $\delta:=\rho_\infty(S)$ exists.

  Then
\[
\lim_{n\to\infty} \sum_{x\in S\cap \sqrt{n}\Omega}
p_n(x/\sqrt{n})=\delta\mathfrak{N}(\Omega).
\]
\end{thm}
\begin{proof}
  The proof is exactly the same as that of Theorem~\ref{thm:vis-free},
  with the only change that instead of Proposition~\ref{prop:zeta} we
  use Proposition~\ref{prop:zeta1}.
\end{proof}

\begin{conv}
  From now and until the end of this section we assume that $k=2$.
\end{conv}

\begin{prop}\label{prop:edensity}
  We have
\[\rho_\infty(U_1(e))=\frac{1}{3}\rho_\infty(U_1)=\frac{1}{3\zeta(2)}.\]
\end{prop}

\begin{proof}
  Let $r,s\ge 1$ be real numbers. For $X,Y\in \{A,O,E\}$ we denote by
  $XY(r,s)$ the number of all $z=(z_1,z_2)\in U_1$ such that $0\le
  z_1< r$, $0\le z_2< s$ and such that the parity of $z_1$ is $X$ and
  the parity of $z_2$ is $Y$. Here $A$ stands for ``any'', $E$ stands
  for ``even'' and $O$ stands for ``odd''.
  
  Let $n\gg 1, m\gg 1$ be integers. Then $AA(n,m)=nm$. We will also
  use $='$ to signify the equality up to an additive error term that
  is $o(nm)$. Note that $EE(n,m)=0$. Then we have
\begin{gather*}
  EO(n,m)= AO(n/2, m)=AA(n/2,m)-AE(n/2,m)=\\
  AA(n/2,m)-OE(n/2,m)=' AA(n/2,m)-EO(n/2,m)='\\
  \frac{1}{2}AA(n,m)-\frac{1}{2}EO(n,m).
\end{gather*}
Therefore
\[
\frac{3}{2} EO(n,m)='\frac{1}{2}AA(n,m)\quad \Rightarrow \quad
EO(n,m)='OE(n,m)='\frac{1}{3}AA(n,m).
\]
Hence $EO(n,m)+OE(n,m)='\frac{2}{3}AA(n,m)$ which implies
\[
OO(n,m)=OO(n,m)+EE(n,m)='\frac{1}{3}AA(n,m).
\]
Since $\rho_\infty(U_1)=\frac{1}{\zeta(2)}$, we have
\[
\lim_{n\to\infty} \frac{AA(n,n)}{n^2}=\frac{1}{\zeta(2)}
\]
Therefore
\[
\lim_{n\to\infty} \frac{EE(n,n)+OO(n,n)}{n^2}=\frac{1}{3\zeta(2)} ,\]
which implies $\rho_\infty(U_1(e))=\frac{1}{3\zeta(2)}$, as required.
\end{proof}

We can now compute the limits for the spherical densities of the set
of visible points for even and odd $n$ tending to infinity for the
case $k=2$.

\begin{thm}\label{thm:spher}
  Let $k=2$. We have
\[
\lim_{m\to\infty} \frac{\gamma_A(2m,
  V_1)}{\gamma_A(2m,F)}=\frac{2}{3\zeta(2)}=\frac{4}{\pi^2}
\]
and
\[
\lim_{m\to\infty} \frac{\gamma_A(2m-1,
  V_1)}{\gamma_A(2m-1,F)}=\frac{8}{\pi^2}.
\]
\end{thm}

\begin{proof}
  
  The proof is essentially the same as that of Theorem~\ref{A'}. We
  present the details for completeness.

  For $c>0$ denote $\Omega_c:=\{x\in \mathbb R^2: ||x||<c\}$. Then
  $\displaystyle\lim_{c\to\infty} \mathfrak{N}(\Omega_c)=1$.  Let
  $\epsilon>0$ be arbitrary. Choose $c>0$ such that
  \[|\mathfrak{N}(\Omega_c)-1|\le \epsilon/3\] and such that
\[
\lim_{n\to\infty}\sum \{ p_n(x): x\in \mathbb Z^2/\sqrt{n} \text{ and
} ||x||\ge c\}\le \epsilon/6.
\]

By Theorem~\ref{thm:vis-free1} and the above formula there is some
$n_0\ge 1$ such that for all $n\ge n_0$ we have
\[
\left| \sum_{x\in \cap \sqrt{n}\Omega_c}
  p_n(x/\sqrt{n})-\frac{1}{3\zeta(2)}\mathfrak{N}(\Omega_c) \right|
\le \epsilon/3
\]
and
\[
\sum \{ p_n(x): x\in \mathbb Z^k/\sqrt{n} \text{ and } ||x||\ge c\}\le
\epsilon/3.
\]

For an \emph{even} $n\ge 2$ let
\begin{gather*}
Q(n):=\\  \frac{\gamma_A(n-1, \{w\in F: \overline{w}\in
U_1(e)\})}{2\gamma_A(n-1,F)}+\frac{\gamma_A(n, \{w\in F:
  \overline{w}\in U_1(e)\})}{2\gamma_A(n,F)}=\\
  \frac{\gamma_A(n, \{w\in F:
  \overline{w}\in U_1(e)\})}{2\gamma_A(n,F)}.
\end{gather*}
In the above equality we use the fact that $n$ is even and all the
points of $U_1(e)$ are even.

For an even $n\ge n_0$ we have
\begin{gather*}
  Q(n)=\\
  \frac{\#\{w\in F: \overline{w}\in U_1(e), |w|_A=n-1 \text{ and }
    ||\overline w||<c\sqrt{n}\}}{2\gamma_A(n-1,F)} +\\ \frac{\#\{w\in
    F: \overline{w}\in U_1(e), |w|_A=n \text{ and } ||\overline
    w||<c\sqrt{n}\}}{2\gamma_A(n,F)}+\\
  \frac{\#\{w\in F: \overline{w}\in U_1(e), |w|_A=n-1 \text{ and }
    ||\overline w||\ge c\sqrt{n}\}}{2\gamma_A(n-1,F)} +\\ 
  \frac{\#\{w\in F: \overline{w}\in U_1(e), |w|_A=n \text{ and }
    ||\overline
    w||\ge c\sqrt{n}\}}{2\gamma_A(n,F)}=\\
  \sum_{x\in U_1(e)\cap \sqrt{n}\Omega_c} p_n(x/\sqrt{n})\ \ +\ 
  \sum_{x\in U_1(e) \cap (\mathbb R^k-\sqrt{n}\Omega_c)}
  p_n(x/\sqrt{n})
\end{gather*}

In the last line of the above equation, the first sum differs from
$\frac{1}{3\zeta(2)}\mathfrak{N}(\Omega_c)$ by at most $\epsilon/3$
since $n\ge n_0$, and the second sum is $\le \epsilon/3$ by the choice
of $c$ and $n_0$. Therefore, again by the choice of $c$, we have
$|Q(n)-\frac{1}{3\zeta(2)}|\le \epsilon$. Since $\epsilon>0$ was
arbitrary, this implies that $\displaystyle\lim_{m\to\infty}
Q(2m)=\frac{1}{3\zeta(2)}$.  Therefore \[ \lim_{m\to\infty}
\frac{\gamma_A(2m, V_1)}{\gamma_A(2m,F)}=2\lim_{m\to\infty}
Q(2m)=\frac{2}{3\zeta(2)}=\frac{4}{\pi^2}.
\]
Together with the conclusion of Theorem~\ref{A} this implies that
\[
\lim_{m\to\infty} \frac{\gamma_A(2m-1,
  V_1)}{\gamma_A(2m-1,F)}=\frac{8}{\pi^2},
\]
as claimed.
\end{proof}

\section{Test elements in the free group of rank two}

A subgroup $H$ of a group $G$ is called a \emph{retract} of $G$ if
there exists a \emph{retraction} from $G$ to $H$, that is, an
endomorphism $\phi:G\to G$ such that $H=\phi(G)$ and that
$\phi|_H=Id_H$. A retract $H\le G$ is \emph{proper} if $H\ne G$ and
$H\ne 1$.

The following result is due to Turner~\cite{Tu}:

\begin{prop}\label{turner}
  Let $F$ be a free group of finite rank $k\ge 2$ and let $w\in F$.
  Then $w$ is a test element in $F$ if and only if $w$ does not belong
  to a proper retract of $F$.
\end{prop}

If $F$ is a free group of rank two, then a proper retract of $F$ is
necessarily cyclic. The following explicit characterization of
retracts in this case is actually Exercise~25 on page~103 of Magnus,
Karrass, Solitar~\cite{MKS}. We present a proof here for completeness.

\begin{lem}\label{mks}
  Let $F$ be a free group of rank two and let $H=\langle h\rangle\le
  F$ be an infinite cyclic subgroup of $F$.

  Then $H$ is a retract of $F$ if and only if there is a free basis
  $\{ a,b \}$ of $F$ such that $h$ can be represented as $h=ac$, where
  $c$ belongs to the normal closure of $b$ in $F$.  In particular, if
  $H$ is a retract of $F$ then $H$ is a maximal cyclic subgroup of
  $F$.
\end{lem}
\begin{proof}
  
  Suppose first that $H$ is a retract of $F$ and that $\phi:F\to F$ is
  a retraction with $\phi(F)=H$. Choose a free basis $x,y$ of $F$.
  Since $H=\langle h\rangle=\langle \phi(x), \phi(y)\rangle\le F$ is
  infinite cyclic, the pair $(x,y)$ is Nielsen equivalent to the pair
  $(h,1)$. Applying the same sequence of Nielsen transformations to
  $(x,y)$ we obtain a free basis $(a,b)$ of $F$ such that $\phi(a)=h$
  and $\phi(b)=1$. Then the kernel of $\phi$ is the normal closure of
  $b$ in $F$.
  
  Since $\phi$ is a retraction onto $H$, we have $\phi(h)=h=\phi(a)$.
  Hence $a^{-1}h\in ker(\phi)$ and therefore $h=ac$, where $c$ belongs
  to the kernel of $\phi$, that is, to the normal closure of $b$, as
  required.
  
  Suppose now that for some free basis $a,b$ of $F$ we have $h=ac$
  where $c$ belongs to the normal closure of $b$ in $F$. Consider the
  endomorphism $\psi:F\to F$ defined by $\psi(a)=h$, $\psi(b)=1$.
  Then, clearly, $\psi(h)=h$ and $\psi$ is a retraction from $F$ to
  $H$.
\end{proof}

We can now obtain an explicit characterization of test elements in
free group $F$ of rank with free basis $A=\{a,b\}$.  We identify the
abelianization of $F$ with $\mathbb Z^2$ so that $\overline a=(1,0)$
and $\overline b=(0,1)$.  If $x\in A$ and $w\in F$, then $w_x$ denotes
the exponent sum on $x$ in $w$ when $w$ is written as a freely reduced
word in $A$ and $\overline w$ denotes the image of $w$ in the
abelianization of $F$. Thus $\overline w=(w_a,w_b)$.

\begin{prop}\label{char-test}
  Let $F$ be a free group of rank two. Let $w\in F$ be a nontrivial
  element that is not a proper power in $F$.
  
  Then $w$ is a test element in $F$ if and only if there exists an
  integer $n\ge 2$ such that $\overline w$ is an $n$-th power in
  $\mathbb{Z}^2$. That is, $w$ is not a test element if and only if
  $w$ is visible in $F$.
\end{prop}

\begin{proof}

  Suppose first that $w$ is a test element but that $\overline w$
  cannot be represented as an $n$-th power in $\mathbb{Z}^2$ for $n\ge
  2$. Then $gcd(w_a,w_b)=1$. Hence there exist integers $p$ and $q$
  such that $pw_a+qw_b=1$.  Consider an endomorphism $\phi:F\to F$
  defined by $\phi(a)=w^p$ and $\phi(b)=w^q$. Then $\phi(w)=w$ and
  $\phi$ is not an automorphism of $F$ since $\phi(F)$ is cyclic.
  Hence, by definition, $w$ is not a test element in $F$, yielding a
  contradiction.

  Suppose now that $\overline w$ is an $n$-th power in $\mathbb{Z}^2$
  for some $n\ge 2$ but that $w$ is not a test element. Then by
  Proposition~\ref{turner} $w$ belongs to an infinite cyclic proper
  retract $H$ of $F$. Since by assumption $w$ is not a proper power in
  $F$, it follows that $w$ generates $H$.  Lemma~\ref{mks} implies
  that for some free basis $(a_1,b_1)$ of $F$ we have $w=a_1c$ where
  $c$ belongs to the normal closure of $b_1$ in $F$. Hence when $w$ is
  expressed as a word in $a_1,b_1$, the exponent sum on $a_1$ in $w$
  is equal to $1$, which contradicts the assumption that $\overline w$
  is an $n$-th power in the abelianization of $F$.
\end{proof}

Note that if $w\in F=F(a,b)$ then $\overline w$ is an $n$-th power in
$\mathbb{Z}^2$ for some $n\ge 2$ if and only if $gcd(w_a,w_b)>1$.  By
convention we set $gcd(0,0)=\infty$.

It is well-known and easy to prove that the set of proper powers in a
free group is negligible~\cite{AO}:

\begin{prop}\label{power} Let $F=F(A)$ be a free group of finite rank
  $k\ge 2$ with free basis $A=\{a_1,\dots,a_k\}$. Let $P$ be the
  set of all nontrivial elements of $F$ that are proper powers.
  
  Then
\begin{gather*}
  \lim_{n\to\infty}
  \frac{\gamma_A(n,P)}{\gamma_A(n,F)}=\lim_{n\to\infty}
  \frac{\rho_A(n,P)}{\rho_A(n,F)}=0.
\end{gather*}
and the convergence in both limits is exponentially fast.
\end{prop}

\begin{proof}[Proof of Theorem~\ref{B}]
  Since $\zeta(2)=\frac{\pi^2}{6}$, Theorem~\ref{B} now follows
  directly from Theorem~\ref{A}, Proposition~\ref{char-test} and
  Proposition~\ref{power}.
    \end{proof}

\section{Open Problems}

As before, let $F=F(A)$ be a free group of rank $k \ge 2$ with free
basis $A=\{a_1,\dots, a_k\}$ and let $\alpha:F\to\mathbb Z^k$ be the
abelianization homomorphism.

\begin{prob}
  Let $k\ge 3$. Is the set of test elements negligible in $F$?  In
  view of Proposition~\ref{turner}, this is equivalent to asking if
  the union of all proper retracts of $F$ is generic in $F$.
\end{prob}

The proof of Proposition~\ref{char-test} shows that a visible element
in $F$ is never a test element and therefore by Theorem~\ref{A} the
asymptotic density of the set of test elements in $F$ is at most
$1-\frac{1}{\zeta(k)}$. For $k\ge 2$ we have
$0<1-\frac{1}{\zeta(k)}<1$ and
$\displaystyle\lim_{k\to\infty}1-\frac{1}{\zeta(k)}=0$.  Thus the
asymptotic density of the set of test elements of $F$ tends to zero as
the rank $k$ of $F$ tends to infinity.

Note that every free factor of $F$ is a retract, but the converse is
not true.  As mentioned in the Introduction, the union of all proper
free factors is negligible in $F$, whereas the union of all proper
retracts is not since every visible element of $F$ belongs to a proper
retract.

\begin{prob}
  For $k\ge 2$ find a subset $S\subseteq \mathbb Z^k$ such that
  $\overline\rho_{||.||_\infty}(S)\ne
  \overline\sigma_A(\alpha^{-1}(S))$.
\end{prob}
Note that if such a set $S$ exists then it is {\it not} invariant
under the action of $SL(k,\mathbb Z)$ in view of Theorem~\ref{A'},

\begin{prob}
  For $w\in F$ define $T(w)=0$ if $\alpha(w)=0$ and define $T(w)$ to
  be the greatest common divisor of the coordinates of $\alpha(w)$ if
  $\alpha(w)\ne 0$. Let $T_n'$ be the expected value of $T$ over the
  sphere of radius $n$ in $F$ with respect to the uniform distribution
  on that sphere and let $T_n=(T_{n-1}'+T_n')/2$.  What can one say
  about the behavior of $T_n$ as $n\to\infty$?
\end{prob}

Using the results of this paper we can show that $\lim_{n\to\infty}
T_n=\infty$ for the case $k=2$. It also seems plausible that for each
$k\ge 3$ we have $\limsup_{n\to\infty} T_n<\infty$ and heuristic
considerations allow us to conjecture that in fact $\lim_{n\to\infty}
T_n=\frac{\zeta(k-1)}{\zeta(k)}$.  A similar question for $\mathbb
Z^2$ has been studied in detail by Diaconis and Erd\"os~\cite{DE}, who
computed the precise asymptotics, as $n\to\infty$, of the expected
value for the greatest common divisor of the coordinates, computed for
the uniform distribution on the $n\times n$-square in $\mathbb Z^2$.

\begin{prob}\label{prob:hyp}
  Let $G$ be a torsion-free \emph{one-ended} word-hyperbolic group. Is
  it true that the set of test elements in $G$ is generic with respect
  to the word-length?
\end{prob}

Although we do not know the asymptotic density of the set of test
elements in a free group of rank $k\ge 3$, one may still expect a
positive answer to Problem~\ref{prob:hyp}, especially if the
hyperbolic group $G$ is not just one-ended but also does not admit
essential $\mathbb Z$-splittings.  In this case the structure of
endomorphisms and automorphisms of $G$ is much more restricted than in
free groups.

As we have seen, the set of proper powers is negligible in free groups
of rank $k\ge 2$ but has positive asymptotic density in free abelian
groups of finite rank. This raises the corresponding question about
free groups in other varieties.  It is possible to show that if $G$ is
a finitely generated nilpotent group and $t\ge 2$ then the set of
$t$-th powers has positive asymptotic density in $G$.

\begin{prob}
  Let $G$ be a finitely generated nonabelian free solvable group.
  What can be said about the asymptotic density of the set of all
  proper powers in $G$?
\end{prob}

\end{document}